\documentclass[12pt,english]{article}

\usepackage[latin9]{inputenc}
\usepackage{geometry}
\usepackage{amsmath}
\usepackage{amssymb}
\usepackage{setspace}
\makeatletter
\@ifundefined{date}{}{\date{}}
\usepackage{array}
\usepackage{mathrsfs}
\usepackage{color}

\makeatother

\usepackage{babel}
\begin{document}

\title{Characterizing the Zeta Distribution via\\
Continuous Mixtures}

\author{Jiansheng Dai,\thanks{WizardQuant Investment Management; email: daiball@yeah.net.}
\enskip{}Ziheng Huang,\thanks{Sixie Capital Management; email: huangziheng1996@126.com.}
\enskip{}Michael R. Powers,\thanks{Corresponding author; 386C Weilun Building, Department of Finance,
School of Economics and Management, and Schwarzman College, Tsinghua
University, Beijing, China 100084; email: powers@sem.tsinghua.edu.cn.} \enskip{}and Jiaxin Xu\thanks{Department of Finance, School of Economics and Management, Tsinghua
University; email: xujx.18@sem.tsinghua.edu.cn.}}

\date{June 4, 2021}
\maketitle
\begin{abstract}
\begin{singlespace}
\noindent We offer two novel characterizations of the Zeta distribution:
first, as tractable continuous mixtures of Negative Binomial distributions
(with fixed shape parameter, $r>0$), and second, as a tractable continuous
mixture of Poisson distributions. In both the Negative Binomial case
for $r\in\left[1,\infty\right)$ and the Poisson case, the resulting
Zeta distributions are identifiable because each mixture can be associated
with a unique mixing distribution. In the Negative Binomial case for
$r\in\left(0,1\right)$, the mixing distributions are quasi-distributions
(for which the quasi-probability density function assumes some negative
values).\medskip{}

\end{singlespace}
\end{abstract}

\section{Introduction}

As part of an investigation of heavy-tailed discrete distributions\footnote{By ``heavy-tailed'', we mean a random variable $X\sim f_{X}\left(x\right)$
for which $E_{X}\left[X^{\alpha}\right]\rightarrow\infty$ for some
$\alpha\in\left(0,\infty\right)$.} in insurance and actuarial science (see Dai, Huang, Powers, and Xu,
2021), the authors derived two new characterizations of the Zeta distribution,
which form the basis for the present article. Within the insurance
context, Zeta random variables sometimes are employed to model heavy-tailed
loss frequencies (i.e., counts of event occurrences, claim submissions,
indemnity payments, etc.; see, e.g., Doray and Arsenault, 2002). However,
their use is even more common in other fields of study, in which they
serve as empirical models for a wide range of discrete processes.
Examples include: the number of word occurrences in a text; various
measures of communication and influence, such as numbers of telephone
calls, emails, and website hits; and physical intensities, such as
numbers of earthquakes and solar flares within specified discrete
categories (see, e.g., Newman, 2005). The Zeta distribution also plays
important roles in analytic number theory, especially with regard
to the distribution of prime numbers and the Riemann Hypothesis (see
Lin and Hu, 2001; and Aoyama and Nakamura, 2012).\footnote{It is important to note that this literature uses the term ``Riemann
Zeta'' to describe random variables formed by taking the negative
of the natural logarithm of conventional ``Zeta'' random variables
(defined on the sample space $x\in\left\{ 1,2,3,\ldots\right\} $).}

The novel formulations of the Zeta distribution presented in this
article are likely to be of interest to researchers in various fields
(as indicated above), especially those exploring social or physical
mechanisms leading to heavy-tailed behavior. The first characterization,
provided in Section 2, shows that Zeta random variables can be expressed
as continuous mixtures of Negative Binomial counts with a fixed shape
parameter, $r>0$. This is accomplished via tractable mixing distributions
that are well behaved for $r\in\left[1,\infty\right)$, but consist
of quasi-distributions (for which the quasi-probability density function
assumes some negative values) for $r\in\left(0,1\right)$. The second
characterization, given in Section 3, converts the mixtures of Negative
Binomial counts from Section 2 into a tractable continuous mixture
of Poisson counts by first expressing each Negative Binomial component
as a continuous mixture of Poisson components. In both the Negative
Binomial case for $r\in\left[1,\infty\right)$ and the Poisson case,
the resulting Zeta distributions are identifiable because the mixing
distributions are unique. In Section 4, we conclude with some final
observations.

\section{Zeta as a Mixture of Negative Binomial Counts}

In this section, we will show how Zeta random variables can be constructed
as continuous mixtures of Negative Binomial counts with a fixed shape
parameter, $r$. To this end, let $X|s\sim\textrm{Zeta}\left(s\right)$
have probability mass function (PMF) $f_{X|s}\left(x\right)=\tfrac{\left(x+1\right)^{-s}}{\zeta\left(s\right)}$
for $x\in\left\{ 0,1,2,\ldots\right\} $ and $s\in\left(1,\infty\right)$,
where $\zeta\left(s\right)={\textstyle \sum_{x=0}^{\infty}}{\displaystyle \left(x+1\right)^{-s}}$
denotes the Riemann zeta function, and let $X\mid r,p\sim\textrm{Negative Binomial}\left(r,p\right)$
have PMF $f_{X\mid r,p}\left(x\mid r,p\right)=\tfrac{\Gamma\left(r+x\right)}{\Gamma\left(r\right)\Gamma\left(x+1\right)}p^{x}\left(1-p\right)^{r}$
for $x\in\left\{ 0,1,2,\ldots\right\} $, fixed $r\in\left(0,\infty\right)$,
and $p\in\left(0,1\right)$.\footnote{The $\textrm{Zeta}\left(s\right)$ and $\textrm{Negative Binomial}\left(r,p\right)$
distributions often are defined on the sample space $x\in\left\{ 1,2,3,\ldots\right\} $
rather than $x\in\left\{ 0,1,2,\ldots\right\} $. However, we have
chosen the latter characterization both because it matches the sample
space of the $\textrm{Poisson}\left(\lambda\right)$ distribution
and because it is the more commonly used formulation in insurance
applications (where it is convenient for loss frequencies to admit
the possibility of $x=0$).} Our basic objective is to identify a mixing random variable, $p|r,s\sim f_{p|r,s}\left(p\right)$
for $p\in\left(0,1\right)$, such that
\begin{equation}
f_{X|s}\left(x\right)=f_{X\mid r,p}\left(x\right)\underset{p}{\wedge}f_{p|r,s}\left(p\right).
\end{equation}

Since the form of the mixing probability density function (PDF), $f_{p|r,s}\left(p\right)$,
differs by the domain of the Negative Binomial shape parameter, $r$,
we must consider the cases of $r=1$, $r\in\left(1,\infty\right)$,
and $r\in\left(0,1\right)$, respectively, in the following three
subsections.

\subsection{The Case of \emph{$\boldsymbol{r=1}$} (Geometric Distribution)}

If the shape parameter equals 1, then the $\textrm{Negative Binomial}\left(r,p\right)$
distribution simplifies to\linebreak{}
$\textrm{Geometric}\left(p\right)$, with PMF $f_{X\mid r=1,p}\left(x\right)=p^{x}\left(1-p\right)$.
Rewriting (1) as
\[
f_{X|s}\left(x\right)=f_{X\mid r=1,p}\left(x\right)\underset{p}{\wedge}f_{p|r=1,s}\left(p\right)
\]
yields
\[
{\displaystyle \int_{0}^{1}}p^{x}\left(1-p\right)f_{p|r=1,s}\left(p\right)dp=\dfrac{\left(x+1\right)^{-s}}{\zeta\left(s\right)}
\]
\[
\Longrightarrow{\displaystyle \int_{0}^{1}}\left(p^{x}-p^{x+1}\right)f_{p|r=1,s}\left(p\right)dp=\dfrac{\left(x+1\right)^{-s}}{\zeta\left(s\right)}
\]
\[
\Longrightarrow E_{p|r=1,s}\left[p^{x+1}\right]=E_{p|r=1,s}\left[p^{x}\right]-\dfrac{\left(x+1\right)^{-s}}{\zeta\left(s\right)},
\]
from which it follows that
\begin{equation}
E_{p|r=1,s}\left[p^{x}\right]=1-\dfrac{1}{\zeta\left(s\right)}{\displaystyle \sum_{i=0}^{x-1}\left(i+1\right)^{-s}}.
\end{equation}

The system (2) then can be used to express the moment-generating function
of $p|r=1,s$ as
\[
E_{p|r=1,s}\left[e^{tp}\right]=1+\dfrac{tE_{p|r=1,s}\left[p\right]}{1!}+\dfrac{t^{2}E_{p|r=1,s}\left[p^{2}\right]}{2!}+\dfrac{t^{3}E_{p|r=1,s}\left[p^{3}\right]}{3!}+\cdots
\]
\begin{equation}
=1+\dfrac{t}{1!}\left(1-\dfrac{1}{\zeta\left(s\right)}\right)+\dfrac{t^{2}}{2!}\left(1-\dfrac{1+2^{-s}}{\zeta\left(s\right)}\right)+\dfrac{t^{3}}{3!}\left(1-\dfrac{{\displaystyle 1+2^{-s}+3^{-s}}}{\zeta\left(s\right)}\right)+\cdots
\end{equation}
\[
=e^{t}-\dfrac{1}{\zeta\left(s\right)}{\displaystyle \sum_{n=1}^{\infty}}\dfrac{t^{n}}{n!}\left(\sum_{i=0}^{n-1}\left(i+1\right)^{-s}\right)
\]
\begin{equation}
=e^{t}-{\displaystyle \sum_{n=1}^{\infty}}\dfrac{t^{n}}{n!}\dfrac{H_{n,s}}{\zeta\left(s\right)},
\end{equation}
where $H_{n,s}={\displaystyle {\textstyle \sum}_{i=0}^{n-1}\left(i+1\right)^{-s}}$
is the $n^{\textrm{th}}$ generalized harmonic number. Alternatively,
(3) may be written as
\[
\dfrac{1}{\zeta\left(s\right)}{\displaystyle \sum_{n=0}^{\infty}}\dfrac{t^{n}}{n!}\left({\displaystyle \sum_{i=n}^{\infty}}{\displaystyle \left(i+1\right)^{-s}}\right)
\]
\begin{equation}
={\displaystyle \sum_{n=0}^{\infty}}\dfrac{t^{n}}{n!}\dfrac{\zeta\left(s,n+1\right)}{\zeta\left(s\right)},
\end{equation}
where $\zeta\left(s,n+1\right)={\displaystyle {\textstyle \sum}_{i=n}^{\infty}\left(i+1\right)^{-s}}$
is the Hurwitz zeta function of order $n+1$.

The series (4) and (5) clearly converge for all $t\in\left(0,\infty\right)$
because both $\tfrac{H_{n,s}}{\zeta\left(s\right)}$ and $\tfrac{\zeta\left(s,n+1\right)}{\zeta\left(s\right)}$
are bounded above as $n\rightarrow\infty$ for all $s\in\left(1,\infty\right)$.
Although these series are not reducible to simpler expressions, it
is possible to write the PDF $f_{p\mid r=1,s}\left(p\right)$ analytically,
as shown in the following result.\medskip{}

\noindent \textbf{Theorem 1:} If $X|s\sim\textrm{Zeta}\left(s\right)$
and $X|r=1,p\sim\textrm{Negative Binomial}\left(r=1,p\right)$, then
there exists a unique mixing random variable, $p|r=1,s$, with PDF
\[
f_{p|r=1,s}\left(p\right)=\dfrac{\left(-\ln\left(p\right)\right)^{s-1}}{\zeta\left(s\right)\Gamma\left(s\right)\left(1-p\right)}
\]
for $p\in\left(0,1\right)$, such that
\[
f_{X|s}\left(x\right)=f_{X\mid r=1,p}\left(x\right)\underset{p}{\wedge}f_{p|r=1,s}\left(p\right).
\]
\textbf{Proof:} First, consider
\[
{\displaystyle \int_{0}^{1}f_{X\mid r=1,p}\left(x\right)f_{p|r=1,s}\left(p\right)dp}={\displaystyle \int_{0}^{1}p^{x}\left(1-p\right)\dfrac{\left(-\ln\left(p\right)\right)^{s-1}}{\zeta\left(s\right)\Gamma\left(s\right)\left(1-p\right)}dp}
\]
\begin{equation}
=\dfrac{1}{\zeta\left(s\right)\Gamma\left(s\right)}{\displaystyle \int_{0}^{1}}p^{x}\left(-\ln\left(p\right)\right)^{s-1}dp.
\end{equation}
Then, using the substitution $y=-\ln\left(p\right)$ in the above
integral, (6) can be rewritten as
\[
\dfrac{1}{\zeta\left(s\right)\Gamma\left(s\right)}{\displaystyle \int_{\infty}^{0}}\left(e^{-y}\right)^{x}y^{s-1}\left(-e^{-y}\right)dy
\]
\[
=\dfrac{1}{\zeta\left(s\right)\Gamma\left(s\right)}{\displaystyle \int_{0}^{\infty}}\left(e^{-y}\right)^{x+1}y^{s-1}dy
\]
\[
=\dfrac{\left(x+1\right)^{-s}}{\zeta\left(s\right)}{\displaystyle \int_{0}^{\infty}}\dfrac{\left(x+1\right)^{s}y^{s-1}e^{-\left(x+1\right)y}}{\Gamma\left(s\right)}dy
\]
\[
=\dfrac{\left(x+1\right)^{-s}}{\zeta\left(s\right)}.
\]
The uniqueness of $f_{p|r=1,s}\left(p\right)$ \textendash{} and equivalently,
the identifiability of $f_{X|s}\left(x\right)$ \textendash{} follows
from the identifiability of Negative Binomial mixtures with fixed
shape parameter $r$ (see Theorem 2.1 of Sapatinas, 1995). $\blacksquare$\medskip{}

One fairly obvious, yet interesting, aspect of the above result is
that it provides a natural connection between two of the simplest
convergent series in mathematical analysis: the infinite geometric
series,
\[
S_{\textrm{G}}=1+\gamma^{-1}+\gamma^{-2}+\gamma^{-3}+\ldots=\dfrac{\gamma}{\gamma-1},
\]
with $\gamma\in\left(1,\infty\right)$; and the zeta function,
\[
S_{\textrm{Z}}=1+2^{-s}+3^{-s}+4^{-s}\ldots=\zeta\left(s\right),
\]
with $s\in\left(1,\infty\right)$. Letting $\tau_{\textrm{G}}\left(n\right)$
and $\tau_{\textrm{Z}}\left(n\right)$ denote the $n^{\textrm{th}}$
terms of these two series, respectively (for $n\in\left\{ 1,2,3,\ldots\right\} $),
and treating $\gamma$ as a random variable defined by the transformation
$\gamma=\tfrac{1}{p}$, it follows from Theorem 1 that
\[
f_{\gamma|s}\left(\gamma\right)=\dfrac{\left(\ln\left(\gamma\right)\right)^{s-1}}{\zeta\left(s\right)\Gamma\left(s\right)\gamma\left(\gamma-1\right)}
\]
 for $\gamma\in\left(1,\infty\right)$ and
\[
\dfrac{\tau_{\textrm{Z}}\left(n\right)}{S_{\textrm{Z}}}={\displaystyle \int_{0}^{1}}p^{n-1}\left(1-p\right)f_{p|r=1,s}\left(p\right)dp
\]
\[
={\displaystyle \int_{0}^{1}}p^{n-1}\left(1-p\right)\dfrac{\left(-\ln\left(p\right)\right)^{s-1}}{\zeta\left(s\right)\Gamma\left(s\right)\left(1-p\right)}dp
\]
\[
=-{\displaystyle \int_{\infty}^{1}\left(\dfrac{1}{\gamma}\right)^{n-1}}\left(1-\dfrac{1}{\gamma}\right)\dfrac{\left(-\ln\left(1/\gamma\right)\right)^{s-1}}{\zeta\left(s\right)\Gamma\left(s\right)\left(1-\dfrac{1}{\gamma}\right)}\gamma^{-2}d\gamma
\]
\[
={\displaystyle \int_{1}^{\infty}}\dfrac{\gamma^{-\left(n-1\right)}}{\dfrac{\gamma}{\gamma-1}}\dfrac{\left(\ln\left(\gamma\right)\right)^{s-1}}{\zeta\left(s\right)\Gamma\left(s\right)\gamma\left(\gamma-1\right)}d\gamma
\]
\[
={\displaystyle \int_{1}^{\infty}}\dfrac{\tau_{\textrm{G}}\left(n\right)}{S_{\textrm{G}}}f_{\gamma|s}\left(\gamma\right)d\gamma
\]
\[
=E_{\gamma|s}\left[\dfrac{\tau_{\textrm{G}}\left(n\right)}{S_{\textrm{G}}}\right].
\]

\subsection{The Case of Fixed \emph{$\boldsymbol{r\in\left(1,\infty\right)}$}}

\noindent When the fixed Negative Binomial shape parameter differs
from 1, the mixing PDF, $f_{p|r,s}\left(p\right)$, generally becomes
more complex, but remains reasonably tractable. The following result
addresses the case of $r\in\left(1,\infty\right)$.\medskip{}

\noindent \textbf{Theorem 2:} If $X|s\sim\textrm{Zeta}\left(s\right)$
and $X|r>1,p\sim\textrm{Negative Binomial}\left(r>1,p\right)$, then
there exists a unique mixing random variable, $p|r>1,s$, with PDF
\[
f_{p|r>1,s}\left(p\right)=\dfrac{r-1}{\zeta\left(s\right)\Gamma\left(s\right)\left(1-p\right)^{r}}{\displaystyle \int_{p}^{1}}\dfrac{\left(\omega-p\right)^{r-2}\left(-\ln\left(\omega\right)\right)^{s-1}}{\omega^{r-1}}d\omega
\]
for $p\in\left(0,1\right)$, such that
\[
f_{X|s}\left(x\right)=f_{X\mid r>1,p}\left(x\right)\underset{p}{\wedge}f_{p|r>1,s}\left(p\right).
\]
\textbf{Proof:} Consider
\[
{\displaystyle \int_{0}^{1}f_{X\mid r>1,p}\left(x\right)f_{p|r>1,s}\left(p\right)dp}
\]
\[
={\displaystyle \int_{0}^{1}\dfrac{\Gamma\left(r+x\right)}{\Gamma\left(r\right)\Gamma\left(x+1\right)}p^{x}\left(1-p\right)^{r}\left[\dfrac{r-1}{\zeta\left(s\right)\Gamma\left(s\right)\left(1-p\right)^{r}}{\displaystyle \int_{p}^{1}}\dfrac{\left(\omega-p\right)^{r-2}\left(-\ln\left(\omega\right)\right)^{s-1}}{\omega^{r-1}}d\omega\right]dp}
\]
\begin{equation}
=\dfrac{1}{\zeta\left(s\right)}\dfrac{\Gamma\left(r+x\right)}{\Gamma\left(r-1\right)\Gamma\left(x+1\right)}\int_{0}^{1}\dfrac{p^{x}}{\Gamma\left(s\right)}\left[{\displaystyle \int_{p}^{1}}\dfrac{\left(\omega-p\right)^{r-2}\left(-\ln\left(\omega\right)\right)^{s-1}}{\omega^{r-1}}d\omega\right]dp.
\end{equation}
Using the substitution $t=-\ln\left(\omega\right)$ in the inside
integral, followed by $y=-\ln\left(p\right)$ in the outside integral,
(7) can be rewritten as
\[
\dfrac{1}{\zeta\left(s\right)}\dfrac{\Gamma\left(r+x\right)}{\Gamma\left(r-1\right)\Gamma\left(x+1\right)}{\displaystyle \int_{\infty}^{0}}\dfrac{\left(e^{-y}\right)^{x}}{\Gamma\left(s\right)}\left[{\displaystyle \int_{y}^{0}}\dfrac{\left(e^{-t}-e^{-y}\right)^{r-2}t^{s-1}}{\left(e^{-t}\right)^{r-1}}\left(-e^{-t}\right)dt\right]\left(-e^{-y}\right)dy
\]
\begin{equation}
=\dfrac{1}{\zeta\left(s\right)}\dfrac{\Gamma\left(r+x\right)}{\Gamma\left(r-1\right)\Gamma\left(x+1\right)}{\displaystyle \int_{0}^{\infty}}{\displaystyle \int_{0}^{y}}\dfrac{1}{\Gamma\left(s\right)}e^{-\left(x+1\right)y}\left(1-e^{t-y}\right)^{r-2}t^{s-1}dtdy.
\end{equation}
Now interchange the order of integration, and let $\xi=y-t$ in the
new inside integral, so (8) becomes
\[
\dfrac{1}{\zeta\left(s\right)}\dfrac{\Gamma\left(r+x\right)}{\Gamma\left(r-1\right)\Gamma\left(x+1\right)}{\displaystyle {\displaystyle \int_{0}^{\infty}}\int_{0}^{\infty}}\dfrac{1}{\Gamma\left(s\right)}e^{-\left(x+1\right)\left(\xi+t\right)}\left(1-e^{-\xi}\right)^{r-2}t^{s-1}d\xi dt
\]
\[
=\dfrac{1}{\zeta\left(s\right)}\dfrac{\Gamma\left(r+x\right)}{\Gamma\left(r-1\right)\Gamma\left(x+1\right)}{\displaystyle {\displaystyle \int_{0}^{\infty}}\dfrac{1}{\Gamma\left(s\right)}t^{s-1}e^{-\left(x+1\right)t}dt\int_{0}^{\infty}}e^{-\left(x+1\right)\xi}\left(1-e^{-\xi}\right)^{r-2}d\xi
\]
\[
=\dfrac{\left(x+1\right)^{-s}}{\zeta\left(s\right)}\dfrac{\Gamma\left(r+x\right)}{\Gamma\left(r-1\right)\Gamma\left(x+1\right)}{\displaystyle {\displaystyle \int_{0}^{\infty}}\dfrac{\left(x+1\right)^{s}}{\Gamma\left(s\right)}t^{s-1}e^{-\left(x+1\right)t}dt\int_{0}^{\infty}}e^{-\left(x+1\right)\xi}\left(1-e^{-\xi}\right)^{r-2}d\xi
\]
\[
=\dfrac{\left(x+1\right)^{-s}}{\zeta\left(s\right)}\dfrac{\Gamma\left(r+x\right)}{\Gamma\left(r-1\right)\Gamma\left(x+1\right)}{\displaystyle \int_{0}^{\infty}}e^{-\left(x+1\right)\xi}\left(1-e^{-\xi}\right)^{r-2}d\xi.
\]
Substituting $q=e^{-\xi}$ into the above integral then yields
\[
\dfrac{\left(x+1\right)^{-s}}{\zeta\left(s\right)}\dfrac{\Gamma\left(r+x\right)}{\Gamma\left(r-1\right)\Gamma\left(x+1\right)}{\displaystyle \int_{1}^{0}}q^{x+1}\left(1-q\right)^{r-2}\left(-\dfrac{1}{q}\right)dq
\]
\[
=\dfrac{\left(x+1\right)^{-s}}{\zeta\left(s\right)}{\displaystyle \int_{0}^{1}}\dfrac{\Gamma\left(x+r\right)}{\Gamma\left(x+1\right)\Gamma\left(r-1\right)}q^{x}\left(1-q\right)^{r-2}dq
\]
\[
=\dfrac{\left(x+1\right)^{-s}}{\zeta\left(s\right)}.
\]
The uniqueness of $f_{p|r>1,s}\left(p\right)$ follows in the same
way as the uniqueness of $f_{p|r=1,s}\left(p\right)$ in the proof
of Theorem 1. $\blacksquare$\medskip{}

It is worth noting that, for the special case of $r=2$, the PDF $f_{p|r>1,s}\left(p\right)$
simplifies to
\[
f_{p|r=2,s}\left(p\right)=\dfrac{1}{\zeta\left(s\right)\Gamma\left(s\right)\left(1-p\right)^{2}}{\displaystyle \int_{p}^{1}}\dfrac{\left(-\ln\left(\omega\right)\right)^{s-1}}{\omega}d\omega
\]
\[
=\dfrac{1}{\zeta\left(s\right)\Gamma\left(s\right)\left(1-p\right)^{2}}\left[\left.-\dfrac{\left(-\ln\left(\omega\right)\right)^{s}}{s}\right|_{p}^{1}\right]
\]
\[
=\dfrac{\left(-\ln\left(p\right)\right)^{s}}{\zeta\left(s\right)\Gamma\left(s+1\right)\left(1-p\right)^{2}},
\]
an analytic form quite similar to $f_{p|r=1,s}\left(p\right)$.

\subsection{The Case of Fixed $\boldsymbol{r\in\left(0,1\right)}$}

For fixed $r\in\left(0,1\right)$, the analysis is similar to that
for $r\in\left(1,\infty\right)$, with one major difference: the mixing
distribution has a quasi-PDF, $f_{p|r<1,s}\left(p\right)$, that assumes
some negative values. We provide the details in the following result.\medskip{}

\noindent \textbf{Theorem 3:} If $X|s\sim\textrm{Zeta}\left(s\right)$
and $X|r<1,p\sim\textrm{Negative Binomial}\left(r<1,p\right)$, then
there exists a mixing quasi-random variable, $p|r<1,s$, with quasi-PDF
\[
f_{p|r<1,s}\left(p\right)=\dfrac{1}{\zeta\left(s\right)\Gamma\left(s\right)\left(1-p\right)^{r}}\left[\left(s-1\right){\displaystyle \int_{p}^{1}}\dfrac{\left(\omega-p\right)^{r-1}\left(-\ln\left(\omega\right)\right)^{s-2}}{\omega^{r}}d\omega\right.
\]
\begin{equation}
\left.+\left(r-1\right){\displaystyle \int_{p}^{1}}\dfrac{\left(\omega-p\right)^{r-1}\left(-\ln\left(\omega\right)\right)^{s-1}}{\omega^{r}}d\omega\right]
\end{equation}
for $p\in\left(0,1\right)$, such that
\[
f_{X|s}\left(x\right)=f_{X\mid r<1,p}\left(x\right)\underset{p}{\wedge}f_{p|r<1,s}\left(p\right).
\]
\textbf{Proof:} After writing
\[
{\displaystyle \int_{0}^{1}f_{X\mid r<1,p}\left(x\right)f_{p|r<1,s}\left(p\right)dp}
\]
\[
={\displaystyle \int_{0}^{1}\dfrac{\Gamma\left(r+x\right)}{\Gamma\left(r\right)\Gamma\left(x+1\right)}p^{x}\left(1-p\right)^{r}\dfrac{1}{\zeta\left(s\right)\Gamma\left(s\right)\left(1-p\right)^{r}}\left[\left(s-1\right){\displaystyle \int_{p}^{1}}\dfrac{\left(\omega-p\right)^{r-1}\left(-\ln\left(\omega\right)\right)^{s-2}}{\omega^{r}}d\omega\right.}
\]
\[
\left.+\left(r-1\right){\displaystyle \int_{p}^{1}}\dfrac{\left(\omega-p\right)^{r-1}\left(-\ln\left(\omega\right)\right)^{s-1}}{\omega^{r}}d\omega\right]dp
\]
\[
=\dfrac{1}{\zeta\left(s\right)}\dfrac{\Gamma\left(r+x\right)}{\Gamma\left(r\right)\Gamma\left(x+1\right)}\int_{0}^{1}p^{x}\left[\dfrac{1}{\Gamma\left(s-1\right)}{\displaystyle \int_{p}^{1}}\dfrac{\left(\omega-p\right)^{r-1}\left(-\ln\left(\omega\right)\right)^{s-2}}{\omega^{r}}d\omega\right.
\]
\begin{equation}
\left.+\dfrac{\left(r-1\right)}{\Gamma\left(s\right)}{\displaystyle \int_{p}^{1}}\dfrac{\left(\omega-p\right)^{r-1}\left(-\ln\left(\omega\right)\right)^{s-1}}{\omega^{r}}d\omega\right]dp,
\end{equation}
the proof that (10) equals $\tfrac{\left(x+1\right)^{-s}}{\zeta\left(s\right)}$
is entirely analogous to the proof that (7) equals $\tfrac{\left(x+1\right)^{-s}}{\zeta\left(s\right)}$
for Theorem 2.

To demonstrate $f_{p|r<1,s}\left(p\right)<0$ for some $p\in\left(0,1\right)$
for all $s\in\left(1,\infty\right)$, it suffices to insert $p=0$
into the right-hand side of (9), revealing $\underset{p\rightarrow0^{+}}{\lim}f_{p|r<1,s}\left(p\right)=-\infty$.
This implies there exists an interval $\left(0,\varepsilon\right)$,
for some $\varepsilon>0$, such that $p\in\left(0,\varepsilon\right)\Longrightarrow f_{p|r<1,s}\left(p\right)<0$.
$\blacksquare$\medskip{}

The reason $f_{p|r<1,s}\left(p\right)<0$ in some neighborhood of
zero is quite intuitive. Essentially, when $r\in\left(0,1\right)$,
the Negative Binomial PMF becomes very steep at the lower end of its
sample space (in the sense that $\underset{r}{\sup}\left(\tfrac{f_{X|r,p}\left(0\right)}{f_{X|r,p}\left(1\right)}\right)=\underset{r\rightarrow0^{+}}{\lim}\left(\tfrac{1}{rp}\right)=\infty$),
and this steepness is aggravated for values of $p$ close to zero.
In this region of the sample space, it is impossible to construct
the much flatter Zeta PMF (for which $\underset{s}{\inf}\left(\tfrac{f_{X|s}\left(0\right)}{f_{X|s}\left(1\right)}\right)=\underset{s\rightarrow1^{+}}{\lim}2^{s}=2$)
as a convex combination of Negative Binomial PMFs. However, if one
can assign negative weight to those Negative Binomial PMFs for which
$p$ is very small, then one can mitigate the impact of small $r$
by offsetting it with negative contributions from small $p$.

\subsection{The Case of Random $\boldsymbol{r}$}

If the Negative Binomial shape parameter, $r$, is not fixed, but
rather is a continuous random variable on $\left(0,\infty\right)$,
then it is possible to express $X|s\sim\textrm{Zeta}\left(s\right)$
as a mixture
\[
f_{X|s}\left(x\right)=f_{X\mid r,p}\left(x\right)\underset{r,p}{\wedge}f_{r,p|s}\left(r,p\right),
\]
for some joint mixing PDF, $f_{r,p|s}\left(r,p\right)$. Unfortunately,
this joint PDF is not unique, and the resulting Zeta distribution
therefore is not identifiable. In fact, if $f_{r,p|s}\left(r,p\right)=f_{p|r,s}\left(p\right)f_{r|s}\left(r\right)$,
where $f_{p|r,s}\left(p\right)$ denotes the mixing PDF given by Theorems
1 or 2 above, then
\[
f_{X|s}\left(x\right)={\displaystyle \int_{0}^{\infty}\int_{0}^{1}f_{X\mid r,p}\left(x\right)f_{r,p|s}\left(r,p\right)drdp}
\]
\[
=\int_{0}^{\infty}\int_{0}^{1}f_{X\mid r,p}\left(x\right)f_{p|r,s}\left(p\right)f_{r|s}\left(r\right)drdp
\]
\[
=\int_{0}^{\infty}\left[\int_{0}^{1}f_{X\mid r,p}\left(x\right)f_{p|r,s}\left(p\right)dp\right]f_{r|s}\left(r\right)dr
\]
\[
=f_{X|s}\left(x\right)\int_{0}^{\infty}f_{r|s}\left(r\right)dr,
\]
so that $f_{r|s}\left(r\right)$ can be any well-defined PDF on $\left(0,\infty\right)$.
Without the identifiability property, it is impossible to probe the
statistical processes generating the mixed distribution of interest
(in our case, $f_{X|s}\left(x\right)$). In particular, one cannot
estimate the parameters of the mixing distribution from observations
of the mixed random variable (see, e.g., Karlis and Xekalaki, 2005).

\section{Zeta as a Mixture of Poisson Counts}

For any choice of fixed $r\in\left(0,\infty\right)$ and $p\in\left(0,1\right)$,
the $\textrm{Negative Binomial}\left(r,p\right)$ random variable
can be expressed as a continuous mixture of $\textrm{Poisson}\left(\lambda\right)$
counts, using a $\textrm{Gamma}\left(r,\beta=\tfrac{1-p}{p}\right)$
mixing distribution with PDF $f_{\lambda|r,\beta=\tfrac{1-p}{p}}\left(\lambda\right)=\left(\tfrac{1-p}{p}\right)^{r}\tfrac{\lambda^{r-1}}{\Gamma\left(r\right)}\exp\left(-\left(\tfrac{1-p}{p}\right)\lambda\right)$
for $\lambda\in\left(0,\infty\right)$, $r\in\left(1,\infty\right)$,
and $\beta=\tfrac{1-p}{p}\in\left(0,\infty\right)$. This mixture
may be written as
\[
f_{X\mid r,p}\left(x\right)=f_{X\mid\lambda}\left(x\right)\underset{\lambda}{\wedge}f_{\lambda|r,\beta=\tfrac{1-p}{p}}\left(\lambda\right).
\]
Substituting the right-hand side of the above equation into (1) yields
\[
f_{X|s}\left(x\right)=\left[f_{X\mid\lambda}\left(x\right)\underset{\lambda}{\wedge}f_{\lambda|r,\beta=\tfrac{1-p}{p}}\left(\lambda\right)\right]\underset{p}{\wedge}f_{p|r,s}\left(p\right),
\]
which, by the associative property of distribution mixing, is equivalent
to
\[
f_{X|s}\left(x\right)=f_{X\mid\lambda}\left(x\right)\underset{\lambda}{\wedge}\left[f_{\lambda|r,\beta=\tfrac{1-p}{p}}\left(\lambda\right)\underset{p}{\wedge}f_{p|r,s}\left(p\right)\right]
\]
\[
=f_{X\mid\lambda}\left(x\right)\underset{\lambda}{\wedge}f_{\lambda|r,s}\left(\lambda\right).
\]
It thus follows that Zeta random variables can be expressed as continuous
mixtures of Poisson counts with the mixing PDF,
\[
f_{\lambda|r,s}\left(\lambda\right)=f_{\lambda|r,\beta=\tfrac{1-p}{p}}\left(\lambda\right)\underset{p}{\wedge}f_{p|r,s}\left(p\right)
\]
\begin{equation}
={\displaystyle \int_{0}^{1}}\dfrac{\left(\dfrac{1-p}{p}\right)^{r}\lambda^{r-1}}{\Gamma\left(r\right)}\exp\left(-\left(\dfrac{1-p}{p}\right)\lambda\right)f_{p|r,s}\left(p\right)dp.
\end{equation}

Although the right-hand side of (11) appears to depend on the Negative
Binomial shape parameter ($r$), that actually is not the case. As
shown by Feller (1943), any mixed-Poisson random variable must possess
the identifiability property, and thus a unique mixing distribution.
Consequently, $f_{\lambda|r,s}\left(\lambda\right)$ must be invariant
over $r$, and may be expressed using the simplest form of the above
integral (i.e., by inserting $r=1$). Therefore,
\[
f_{\lambda|r,s}\left(\lambda\right)=f_{\lambda|s}\left(\lambda\right)
\]
\[
={\displaystyle \int_{0}^{1}}\left(\dfrac{1-p}{p}\right)\exp\left(-\left(\dfrac{1-p}{p}\right)\lambda\right)\dfrac{\left(-\ln\left(p\right)\right)^{s-1}}{\zeta\left(s\right)\Gamma\left(s\right)\left(1-p\right)}dp
\]
\begin{equation}
=\dfrac{1}{\zeta\left(s\right)\Gamma\left(s\right)}{\displaystyle \int_{0}^{1}}\dfrac{1}{p}\left(-\ln\left(p\right)\right)^{s-1}\exp\left(-\left(\dfrac{1-p}{p}\right)\lambda\right)dp.
\end{equation}
Applying the substitution $y=\tfrac{1-p}{p}$ to the above integral,
(12) can be rewritten as
\[
\dfrac{1}{\zeta\left(s\right)\Gamma\left(s\right)}{\displaystyle \int_{\infty}^{0}}\left(y+1\right)\left(\ln\left(y+1\right)\right)^{s-1}\exp\left(-\lambda y\right)\left[\dfrac{-1}{\left(y+1\right)^{2}}\right]dy
\]
\[
=\dfrac{1}{\zeta\left(s\right)\Gamma\left(s\right)}{\displaystyle \int_{0}^{\infty}}\dfrac{1}{y+1}\left(\ln\left(y+1\right)\right)^{s-1}\exp\left(-\lambda y\right)dy,
\]
which is not further reducible.

\section{Conclusions}

The present article provided two novel characterizations of $X|s\sim\textrm{Zeta}\left(s\right)$.
First, we showed that these random variables can be expressed as tractable
continuous mixtures of $X|r,p\sim$ \linebreak{}
$\textrm{Negative Binomial}\left(r,p\right)$ with fixed shape parameter
(\emph{$r$}); that is,
\[
f_{X|s}\left(x\right)=f_{X\mid r,p}\left(x\right)\underset{p}{\wedge}f_{p|r,s}\left(p\right),
\]
where
\[
f_{p|r,s}\left(p\right)=\dfrac{1}{\zeta\left(s\right)\Gamma\left(s\right)\left(1-p\right)^{r}}\times\begin{cases}
\left(s-1\right){\displaystyle \int_{p}^{1}}\dfrac{\left(\omega-p\right)^{r-1}\left(-\ln\left(\omega\right)\right)^{s-2}}{\omega^{r}}d\omega\\
\quad+\left(r-1\right){\displaystyle \int_{p}^{1}}\dfrac{\left(\omega-p\right)^{r-1}\left(-\ln\left(\omega\right)\right)^{s-1}}{\omega^{r}}d\omega & \textrm{for }r\in\left(0,1\right)\\
\left(-\ln\left(p\right)\right)^{s-1} & \textrm{for }r=1\\
\left(r-1\right){\displaystyle \int_{p}^{1}}\dfrac{\left(\omega-p\right)^{r-2}\left(-\ln\left(\omega\right)\right)^{s-1}}{\omega^{r-1}}d\omega & \textrm{for }r\in\left(1,\infty\right)
\end{cases},
\]
$f_{p|r,s}\left(p\right)=f_{p|r=1,s}\left(p\right)$ and $f_{p|r,s}\left(p\right)=f_{p|r>1,s}\left(p\right)$
are unique PDFs, and $f_{p|r,s}\left(p\right)=f_{p|r<1,s}\left(p\right)$
is a quasi-PDF (i.e., with some negative values). Next, based on the
fact that $\textrm{Negative Binomial}\left(r,p\right)$ random variables
can be constructed as mixtures of $X|\lambda\sim\textrm{Poisson}\left(\lambda\right)$,
with a $\textrm{Gamma}\left(r,\beta=\tfrac{1-p}{p}\right)$ mixing
distribution, we showed that Zeta random variables also can be expressed
as a unique and tractable continuous mixture of Poisson counts; that
is,
\[
f_{X|s}\left(x\right)=f_{X\mid\lambda}\left(x\right)\underset{\lambda}{\wedge}f_{\lambda|s}\left(\lambda\right),
\]
where
\[
f_{\lambda|s}\left(\lambda\right)=\dfrac{1}{\zeta\left(s\right)\Gamma\left(s\right)}{\displaystyle \int_{0}^{\infty}}\dfrac{1}{y+1}\left(\ln\left(y+1\right)\right)^{s-1}\exp\left(-\lambda y\right)dy.
\]

The appearance of quasi-PDFs in the Negative Binomial case for $r\in\left(0,1\right)$
was somewhat unexpected, but \textendash{} as argued in Subsection
2.3 \textendash{} has a fairly intuitive explanation. Therefore, it
is natural to consider whether or not other heavy-tailed discrete
random variables formed as mixtures of Negative Binomial counts also
involve quasi-distributions. One obvious family to consider is $X|b\sim\textrm{Yule}\left(b\right)$,
with PMF $f_{X|b}\left(x\right)=\tfrac{b\Gamma\left(b+1\right)\Gamma\left(x+1\right)}{\Gamma\left(x+b+2\right)}$
for $x\in\left\{ 0,1,2,\ldots\right\} $ and $b\in\left(0,\infty\right)$,
which approximates $X|s\sim\textrm{Zeta}\left(s\right)$ for $b=s-1$.
In fact, $\textrm{Yule}\left(b\right)$ random variables may be expressed
as mixtures of $\textrm{Negative Binomial}\left(r=1,p\right)$ (i.e.,
$\textrm{Geometric}\left(p\right)$) counts, using $f_{p|a=1,b}\left(p\right)\sim\textrm{Beta}\left(a=1,b\right)$
as the mixing distribution; that is,
\[
f_{X|b}\left(x\right)=f_{X\mid r=1,p}\left(x\right)\underset{p}{\wedge}f_{p|a=1,b}\left(p\right).
\]

For the Yule distribution, $\underset{b}{\inf}\left(\tfrac{f_{X|b}\left(0\right)}{f_{X|b}\left(1\right)}\right)=\underset{b\rightarrow0^{+}}{\lim}\left(b+2\right)=2$,
which is identical to the corresponding result for the Zeta distribution
presented at the end of Subsection 2.3. Consequently, one might anticipate
that constructing Yule random variables as mixtures of Negative Binomial
counts with $r\in\left(0,1\right)$ would require quasi-distributions,
as in the Zeta case. In Dai, Huang, Powers, and Xu (2021), we not
only show that this is indeed true, but also provide complete results
for the Yule distribution analogous to those of Theorems 1-3 and Section
3 of the present article.

\end{document}